\documentclass[11pt]{amsart}
\usepackage{amsmath,amssymb,amsthm,amscd}

\def\p{\partial}
\def\R{\mathbb{R}}

\def\l{\lambda}

\numberwithin{equation}{section}

\newtheorem{prop}{Proposition}[section]
\newtheorem{theo}[prop]{Theorem}

\newtheorem{conj}[prop]{Conjecture}
\newtheorem{pl}[prop]{Problem}

\begin{document}
\title{Local solution and extension to the Calabi flow}
\author{Weiyong He}
\address{Department of Mathematics\\
University of British Columbia\\
Vancouver, B.C., V6T 1Z2\\
Canada}
\email{whe@math.ubc.ca}
\thanks{The author is
partially supported by a PIMS postdoctoral fellowship.}
\date{\today}
\maketitle

\section{Introduction}
Let $M$ be a compact K\"ahler manifold of complex dimension $n$ and  let $g$ be a smooth K\"ahler metric on $M$. A K\"ahler metric $g$ can be written as a Hermitian matrix valued function  in  holomorphic coordinates $(z^1, \cdots, z^n)$, \[ g=g_{i\bar j}dz^i\otimes d\bar z^j.
\]
The corresponding K\"ahler form $\omega=\sqrt{-1}g_{i\bar j}dz^i\wedge d\bar z^{j}$ is a closed $(1, 1)$ form. So a K\"ahler class $[\omega]$ defines a nontrivial cohomology class. For simplicity, we use both $g$ and $\omega$ to denote the same K\"ahler metric. 
For any K\"ahler metric $\omega_\varphi$ in the same K\"ahler class $[\omega]$, we can write 
\[
\omega_\varphi=\omega+\sqrt{-1} \partial \bar \partial \varphi,
\]
where $\varphi$ is a real valued function and it is called a K\"ahler potential.
The Calabi flow is defined by
\begin{equation}\label{1-1}
\begin{split}
\frac{\partial \varphi}{\partial t}&=R_\varphi-\underline{R},\\
\varphi(0)&=\varphi_0,
\end{split}
\end{equation}
where $R_\varphi$ is the scalar curvature of $\omega_\varphi,$ $\underline{R}$ is the average of the scalar curvature and it is a constant depending only on $(M, [\omega]).$ In local coordinates, the Ricci curvature of a K\"ahler metric $g$ is given by
\[
R_{i\bar j}=-\frac{\partial^2 }{\partial z^i\partial \bar z^{j}}\log \det (g_{k\bar l}),
\]
and the scalar curvature is given by
\begin{equation}\label{E-1-2}
R=-g^{i\bar j} \frac{\partial^2 }{\partial z^i\partial \bar z^{j}}\log \det (g_{k\bar l}).
\end{equation}
Equation \eqref{1-1} is the gradient flow of the Calabi energy,   aiming to find an extremal metric or a constant scalar curvature metric on a compact K\"ahler manifold within  a fixed K\"ahler class \cite{Calabi82}, which is an important problem in K\"ahler geometry. 
In the case of Riemann surfaces ($n=1$), P. Chrusci\'el \cite{chru} showed that the Calabi flow exists for all time and converges to a constant Gaussian curvature metric, making use of the existence of such a metric and  the Bondi mass in general relativity. Without using the Bondi mass, Chen \cite{chen01}, Struwe \cite{Struwe} reproved Chrusci\'el's theorem. Recently Chen-Zhu \cite{chenzhu07} (for two sphere) removed the assumption that there is a constant Gaussian curvature metric.  While in higher dimensions ($n\geq 2$), the nonlinearity of the Calabi flow becomes more acute and the analytic difficulty becomes more dauting. But recently some progress has been made \cite{Chenhe} \cite{tw} \cite{szekelyhidi} \cite{Chenhe1} \cite{Fine} etc...\\

From \eqref{E-1-2}, we can see that (\ref{1-1}) is a fourth order quasi-linear parabolic equation.  A long standing problem about the Calabi flow ($n\geq 2$) is whether it exists for all  time or not. X. Chen conjectured that the Calabi flow always exists for all time on a compact K\"ahler manifold with smooth initial data. 

\begin{conj}[X. Chen] The Calabi flow exists for all time for any smooth K\"ahler metric on a compact K\"ahler manifold $(M, [\omega])$.
\end{conj}

The motivation of this conjecture relies in part on the fact that the Calabi flow is a distance contracting flow over the space of K\"ahler metrics \cite{Calabi-Chen}, \cite{chen03}. X. Chen then suggested to seek the optimal condition on initial data for the local solution to the Calabi flow. In particular,  the distance over the space of K\"ahler metrics makes sense for  $L^\infty$ K\"ahler metrics. It motivates  to consider the local solution to the Calabi flow for $L^\infty$ initial metrics. X. Chen asked
\begin{pl}[X. Chen]Can  the Calabi flow start with a $L^\infty$ K\"ahler metric?
\end{pl}

In a joitn work with X. Chen \cite{Chenhe}, by using standard  quasilinear parabolic method,  we proved that the short time solution to the Calabi flow exists for initial potential $\varphi_0\in
C^{3, \alpha}(M)$, namely when the initial metric \[g_{\varphi_0}=g+\frac{\partial^2 \varphi_0}{\partial z^i\partial \bar z^{j}}dz^i\otimes d\bar z^j\] is in $C^{1, \alpha}$. 
The quasilinear parabolic nature of the equation allows the solution with a sigularity at time $t=0$ (the curvature is not  pointwise well defined for $\varphi_0\in C^{3, \alpha}$). 
In this note we consider the initial value in $C^{2, \alpha}$ (the initial metric is in $C^\alpha$).  
We prove that for any smooth metric $g\in (M, [\omega])$, there exists a small $C^\alpha$ neighborhood around $g$ ($0<\alpha<1$), such that the short time solution to the Calabi flow exists for all $C^\alpha$ metric in this neighborhood (see Theorem \ref{T-1-1}).  As an application, we can prove that if $g$ itself is a constant scalar curvature metric, there exists a small $C^\alpha$ neighborhood of $g$, such that for any initial metric in this neighborhood, the Calabi flow exists for all time and converges to a constant scalar curvature metric in the same class. This result is proved previously in \cite{Chenhe} for $C^{1, \alpha}$ small neighborhood. 

In a joint work with X. Chen \cite{Chenhe}, we proved that the flow can be extended once the Ricci curvature is uniformly bounded along the Calabi flow. 
In this direction we can prove that on a K\"ahler surface, the Calabi flow can be extended once the metrics are bounded  in $L^\infty$ along the Calabi flow (see Theorem \ref{T-2}), by using a blowing up argument as in \cite{Chenhe1}. In particular, it implies that if we can obtain the $L^\infty$ estimate of the metric $g$ (or the $C^{1, 1}$ bound of the potential $\phi$), we can get that all higher derivatives of $g$ (or the higher derivatives of $\phi$ respectively) are bounded. We can recall the seminal Evans-Krylov theory for second order fully nonlinear equations (elliptic or parabolic type). The theory asserts that when the (elliptic) operator is  concave, then one can obtain the H\"older estimates of the second order derivatives, and then higher regularity follows from the standard boot-strapping argument.  Our result can be viewed as obtaining higher order estimates  only assuming $C^{1, 1}$ bound of K\"ahler potential $\varphi$, by a geometric method for the Calabi flow, a fourth order nonlinear parabolic equation, on K\"ahler surfaces.
 
\section{Local solution for $C^\alpha$ initial metrics}
To state and prove our theorem,  we need to define some Banach spaces. 
The setting is similar to that in \cite{CS}, \cite{Chenhe}.
In the following we assume
$\theta \in (0,1]$, and $E$ is a Banach space, $J=[0, T]$ for some
$T>0$. We consider functions defined on $\widetilde{J}=(0, T]$
with a prescribed singularity at $0$. Set
\begin{eqnarray*} 
C_{1-\theta}(J, E):=\left\{u\in C(\widetilde{J},
E); [t\rightarrow t^{1-\theta}u] \in C(J, E), \lim_{t\rightarrow
0^+}t^{1-\theta}|u(t)|_E=0\right\},\end{eqnarray*}

\begin{eqnarray*}|u|_{C_{1-\theta}(J, E)}:=\sup_{t\in
\widetilde{J}}t^{1-\theta}|u|_E,\end{eqnarray*} 
and
\begin{eqnarray*} C^1_{1-\theta}(J, E):=\left\{u \in
C^1(\widetilde{J}, E); u, \dot{u} \in C_{1-\theta}(J,
E)\right\}.\end{eqnarray*} 

One can verify that $C_{1-\theta}(J,
E)$, equipped with the norm $\|.\|_{C_{1-\theta}(J, E)}$, is a
Banach space and $C^1_{1-\theta}(J, E)$ is a subspace.  We also
set
\[
C_0^1(J, E)=C^1(J, E),~~~C_0(J, E)=C(J, E).
\]
Given $x\in E$, we use the notation $B_{E}(x, r)$ to denote  the
ball in the Banach space $E$ with the center $x$ and the
radius $r$.

Let $E_1, E_0$ be two Banach spaces such that $E_1$ is
continuously, densely embedded into $E_0$. Set 
\begin{eqnarray*}
E_0(J):=C_{1-\theta}(J, E_0),
\end{eqnarray*} 
and
\begin{eqnarray*} 
E_1(J):=C^1_{1-\theta}(J, E_0)\cap
C_{1-\theta}(J, E_1),
\end{eqnarray*}
 where
\begin{eqnarray*}
|u|_{E_1(J)}:=\sup_{t\in
\widetilde{J}}t^{1-\theta}\left(|\dot{u}|_{E_0}+|u|_{E_1}\right).
\end{eqnarray*}
We also use the notation 
\[
E_{\mu}:=(E_0, E_1)_{\mu}, \mu\in (0, 1),
\]
for the continuous interpolation spaces of DaPrato and Grisvard \cite{DP1}.

Recall that if $k \in \mathbb{N}, \alpha \in (0, 1)$, the H\"older
space $C^{k, \alpha}$ is the Banach space of all $C^k$ functions
$f: \mathbb{R}^n\rightarrow \mathbb{R}$ which have finite H\"older
norm. The subspace of $C^{\infty}$ function in $C^{k, \alpha}$ is
not dense under the H\"older norm. One defines the little H\"older
space $c^{k, \alpha}$ to be the closure of smooth functions in the
usual H\"older space $C^{k,\alpha}$. And one can verify $c^{k,
\alpha}$ is a Banach space and that $c^{l,\beta}\hookrightarrow
c^{k, \alpha}$ is a continuous and dense imbedding for $k\leq l$
and $0<\alpha<\beta<1$. These definitions can be extended to
functions on a smooth manifold $M$ naturally \cite{BJ}. For our
purpose, the key fact \cite{Tr} about the continuous interpolation
spaces is that for $k\leq l$ and $0<\alpha<\beta<1$, and
$0<\mu<1$, there is a Banach space isomorphism
\begin{equation}\label{1-2}(c^{k, \alpha}, c^{l, \beta})_{\mu}\cong
c^{\mu l+(1-\mu)k+\mu
\beta+(1-\mu)\alpha},
\end{equation} 
provided that the exponent
$\mu l+(1-\mu)k+\mu \beta+(1-\mu)\alpha$ is not an
integer.

From now on we set $E_0=c^{\alpha}(M, g), E_1=c^{4, \alpha}(M, g)$, and take $\theta=\mu=1/2$, then $E_{1/2}=c^{2, \alpha}(M, g),$ where $g$ is a fixed K\"ahler metric.
For simplicity we will use $c^{\alpha}(M)$ etc if there is no confusion.   One can write
\[
R_\varphi=-A(\nabla^2\varphi)\varphi+f(\nabla \varphi,
\nabla^2\varphi, \nabla^3\varphi),
\]
where the operator $A$ is given by
\[
A(\nabla^2\varphi)w=g^{i\bar j}_{\varphi}g^{k\bar
l}_{\varphi}\nabla_{i}\nabla_{\bar j}\nabla_{k}\nabla_{\bar l}w.
\]
Note that $A(\nabla^2\varphi)$ involves only the
second derivative of $\varphi$.  This fact is important for the initial potential in $c^{2, \alpha}$.  As a quasi-linear equation,  the equation (\ref{1-1}) has lower order term  $f(\nabla\varphi, \nabla^2\varphi,
\nabla^3\varphi)$ which involves the third derivatives of $\varphi$. The
standard quasi-linear theory can apply for $\varphi_0\in
c^{3, \alpha}.$ However, computation shows that $f(\nabla\varphi, \nabla^2\varphi, \nabla^3\varphi)$ involves only quadratic terms of $\nabla^3\varphi$, which allows us to obtain
\begin{theo}\label{T-1-1}For any smooth K\"ahler metric $g$ on $M$, 
there exist positive constants $T=T(g),
\varepsilon=\varepsilon(g)$ and $c=c(g)$ such that for any initial value $x\in E_{1/2}= c^{2, \alpha}$, $|x|_{c^{2, \alpha}}\leq \varepsilon$
the Calabi flow equation (\ref{1-1}) has a unique solution
\[
\varphi(\cdot, x)\in C^1_{1/2}([0, T], E_0)\cap C_{1/2}([0,
T], E_1)
\]
on $[0, T].$  Moreover $\varphi\in C([0, T], E_{1/2})$ and for $x, y\in B_{E_{1/2}}(0, \varepsilon)$, we have 
\[
\begin{split}
\|\varphi(\cdot, x)-\varphi(\cdot, y)\|_{C([0, T],
E_{1/2})}&\leq c\|x-y\|_{E_{1/2}} \\
\|\varphi(\cdot, x)-\varphi(\cdot, y)\|_{E_1([0, T])}&\leq
c\|x-y\|_{E_{1/2}}.
\end{split}
\]
\end{theo}

As an application, we can improve the stability theorem in \cite{Chenhe} as follows,
\begin{theo}\label{T-3}
Suppose $g$ is a constant scalar curvature  metric in $[\omega]$ on $M$.
If the initial metric $g_{\varphi_0}=g+\partial \bar \partial \varphi_0$ satisfies $|\varphi_0|_{c^{2, \alpha}(M)}<\varepsilon$, where $\varepsilon=\varepsilon (g)$ is small enough, then the Calabi flow  exists for all time and
$g(t)=g_{\varphi(t)}$ converges to a constant scalr curvature  metric $g_\infty$  exponentially fast in the
same class $[\omega]$ in $C^{\infty}$ sense. 
\end{theo}
Let $x\in B_{E_{1/2}}(0, \varepsilon)$ for a fixed small  $\varepsilon>0$.   Recall \[ E_1(J)=C^1_{1/2}(J, E_0)\cap C_{1/2}(J,
E_1).
\]
For any $x\in B_{E_{1/2}}(0, \varepsilon)$, set
\[
V_x(J)=\left\{v\in E_1(J);\; v(0)=x, \|v\|_{C(J, E_{1/2})}\leq
\varepsilon_0\right\}\cap B_{E_1(J)}(0, \varepsilon_0)
\]
and equip this set with the topology of $E_1(J)$, where
$\varepsilon_0=C_0\varepsilon$ and $C_0$ is a fixed constant depending only on $g$ and will be determined later. We will show that $V_x(J)$ is not empty if $T$ is small enough. 
Let $v\in V_x(J)$ be given.
Consider the following linear parabolic equation for $u$ such that
for $t\in \tilde J=(0, T]$,
\begin{equation}\label{E-2-2}
\begin{split}
\frac{\partial u}{\partial t}+\triangle^2_gu&= f(v),\\
u(0)&=x,\end{split}
\end{equation}
where
\[
f(v)=\triangle_g^2v+R_v-\underline{R}.
\]
We first show that the linear
equation (\ref{E-2-2}) has a unique solution $u\in V_x(J)$ for any
$v\in V_x(J)$ provided $T, \varepsilon$ sufficiently small. So we can define a map
\[\Pi_x: v\rightarrow u\] from $V_x(J)$ to $V_x(J)$. Note $V_x(J)$
is a close subset of the Banach space $E_1(J)$. Then we show that
the map $\Pi_x$ is a contracting map when $T, \varepsilon$ are both small. It follows that  $\Pi_x$ has a fixed point in $V_x(J)$ by contract mapping theorem and the fixed point is the desired solution.  The smoothing property will
follow simultaneously. We will need some facts of the linear parabolic theory, in particular the analytic semigroup generated by $\triangle_g^2$. The results hold for more general linear operators. One can find a nice reference such as in \cite{CS}, \cite{Lunardi}. 

Consider the linear equation for some $x\in E_0$
\begin{equation}\label{linear}
\begin{split}
\frac{\partial u}{\partial t}&=-\triangle_g^2u,\\
u(0)&=x.
\end{split}
\end{equation}
Then $\triangle_g^2: E_1\rightarrow E_0$ defines an analytic semigroup $e^{-t\triangle^2_g}: E_0\rightarrow E_1$. The solution of \eqref{linear} is given by
$u(t)=e^{-t\triangle_g^2}x$ for $t\in J.$ For $x\in E_{1/2}$, we have the following equivalent norm
\begin{equation}\label{E-5}
|x|_{E_{1/2}}:=\sup_{s>0}s^{1/2}\left|\triangle_g^2\left(e^{-t\triangle_g^2}x\right)\right|_{E_0}.
\end{equation}
Moreover $[ t\longmapsto e^{-t\triangle_g^2}x]\in E_1(J)$ and there exists a constant $c_1>0$ independent of $J$ such that for any $t\in J$
\begin{equation}\label{E-6}
\left|e^{-t\triangle_g^2}x\right|_{E_1(J)}\leq c_1|x|_{E_{1/2}}.
\end{equation}
 If $u\in E_1(J)$ with $u(0)=0$, then there exists a constant $c_2$ independent of $J$ such that
\begin{equation}\label{E-7}
|u|_{C_1(J, E_{1/2})}\leq c_2|u|_{E_1(J)}. 
\end{equation}

The inequalities (\ref{E-5})---(\ref{E-7}) hold for more general settings and one can find the proof in \cite{CS} (see Remark 2.1, Lemma 2.2 in \cite{CS}). Also if $T$ is small, by the strong continuity of the semigroup $\{e^{-t\triangle_g^2}, t\geq 0\}$, we can get that
\begin{equation}\label{E-8}
\left|e^{-t\triangle^2_g}x-x\right|_{E_{1/2}}\leq \frac{1}{4}\varepsilon_0.
\end{equation}
In particular, for any $x\in B_{E_{1/2}}(0, \varepsilon)$,  by (\ref{E-6}),
\[
\left|e^{-t\triangle_g^2}x\right|_{E_1(J)}\leq c_1|x|_{E_{1/2}}.
\]
We  choose $C_0= 2c_1+1$ , then we have

\[\left|e^{-t\triangle_g^2}x\right|_{E_1(J)}\leq \varepsilon_0/2 .\] 
Also by (\ref{E-8}) we get 
\[\left|e^{-t\triangle_g^2}x\right|_{E_{1/2}}\leq |x|_{E_{1/2}}+\left|e^{-t\triangle_g^2}x-x\right|_{E_{1/2}}<\varepsilon_0\] provided $T$ is small.  It follows that $V_x(J)$ is not empty.

Now we are in the position to prove Theorem \ref{T-1-1}. 
\begin{proof}For any $v\in V_x(J)$, the linear parabolic equation (\ref{E-2-2}) has the
unique solution which takes the formula
\[u(t)=\left(e^{-t\triangle^2_g}\right)x+\left(Kf\right)(t),\]
where
\[\left(Kf\right)(t)=\int_0^te^{-(t-s)\triangle_g^2}f(v(s))ds.
\]
It is clear that $K: E_0(J)\rightarrow E_1(J)$ is a linear
operator. We can define the norm of $K$ by $\|K\|$.  We can get for some $c_3$ depending only on $g$,
\begin{equation}\label{E-2-3}|Kf|_{E_1(J)}\leq \|K\| |f|_{E_0(J)}\leq
c_3|f|_{E_0(J)}.\end{equation}

First we need to show $u\in  V_x(J)$. 
We compute
\begin{equation}\label{E-2-5}
\begin{split}
f(v)=&\triangle_g^2v+R_v-\underline{R}\\
=&-{g^{i\bar j}_v}\p_i\p _{\bar j}(\log {\det (g_{k\bar l}+v_{k\bar l})})+g^{i\bar j}\p_i\p_{\bar j} (g^{k\bar l}\p_k\p_{\bar l} v)-\underline{R}\\
=&g^{i\bar j}_v g^{k\bar q}_v g^{p\bar l}_v(\p_i g_{p\bar q}+\p_i\p_p\p_{\bar q} v)(\p_{\bar j} g_{k\bar l}+\p_j\p_k\p_{\bar l} v)-g^{i\bar j}_vg^{k\bar l}_v \p_i\p_{\bar j}\p_k\p_{\bar l}v\\
&+g^{i\bar j}g^{k\bar l}\p_i\p_{\bar j}\p_k\p_{\bar l} v-g^{i\bar j} \p_ig^{k\bar l}\p_{\bar j} v_{k\bar l}-g^{i\bar j} g^{i\bar j}\p_{\bar j}g^{k\bar l}\p_i v_{k\bar l}\\
&+g^{i\bar j}\p_i\p_{\bar j} (g^{k\bar l}) v_{k\bar l}-\underline{R}.
\end{split}
\end{equation}

For simplicity, from now on we use $C$ to denote a uniformly bounded constant independent of $\varepsilon, T$ and can vary line by line. 
We can also compute that 
\begin{equation}\label{E-2-10}
\begin{split}
g^{i\bar j}g^{k\bar l}\p_i\p_{\bar j}\p_k\p_{\bar l} v-g^{i\bar j}_vg^{k\bar l}_v \p_i\p_{\bar j}\p_k\p_{\bar l}v=&\left((g^{i\bar j}-g^{i\bar j}_v)g^{k\bar l} +g^{i\bar j}_v(g^{k\bar l}-g^{k\bar l}_v)\right) \p_i\p_{\bar j}\p_k\p_{\bar l} v\\
=& \left(g^{k\bar l}g^{i\bar q}g^{p\bar j}_v+g^{i\bar j}_vg^{k\bar q} g^{p\bar l}_v\right)\p_p\p_{\bar q}v\p_i\p_{\bar j}\p_k\p_{\bar l} v.
\end{split}
\end{equation}

It follow from \eqref{E-2-5} and \eqref{E-2-10} that

\begin{equation}\label{E-2-6}
\begin{split}
f(v)=&\left(g^{k\bar l}g^{i\bar q}g^{p\bar j}_v+g^{i\bar j}_vg^{k\bar q} g^{p\bar l}_v\right)\p_p\p_{\bar q}v\p_i\p_{\bar j}\p_k\p_{\bar l} v-g^{i\bar j} \p_ig^{k\bar l}\p_{\bar j} v_{k\bar l}-g^{i\bar j} g^{i\bar j}\p_{\bar j}g^{k\bar l}\p_i v_{k\bar l}\\
&+g^{i\bar j}_v g^{k\bar q}_v g^{p\bar l}_v(\p_i g_{p\bar q}+\p_i\p_p\p_{\bar q} v)(\p_{\bar j} g_{k\bar l}+\p_j\p_k\p_{\bar l} v)+g^{i\bar j}\p_i\p_{\bar j} (g^{k\bar l}) v_{k\bar l}-\underline{R}.
\end{split}
\end{equation}

By (\ref{E-2-6}), we can get that
\begin{eqnarray*} \sqrt{t}|f(v(t))|_{E_0}&\leq & C\sqrt{t}(|v|^2_{c^{3, \alpha}}+|v|_{c^{2,\alpha} }|v|_{c^{4, \alpha}}+1)\\
&\leq & C \sqrt{t} (|v|_{c^{2,\alpha}} |v|_{c^{4, \alpha}}+1)\\
&\leq& C\varepsilon_0^2+C\sqrt{T}, ~~~t\in (0, T],
\end{eqnarray*}
where we have used the convexity of H\"older norms,
\[
|v|^2_{c^{3, \alpha}}\leq C |v|_{c^{2,\alpha} }|v|_{c^{4, \alpha}}.
\]
It follows that $f(v)\in E_0(J).$ By (\ref{E-8}) and (\ref{E-2-3}), we can get that
\begin{eqnarray}\label{E-2-7}
|u|_{C(J, E_{1/2})}&\leq&
\left|\left(e^{-t\triangle_g^2}\right)x\right|_{C(J,
E_{1/2})}+\left|(Kf)(t)\right|_{C(J, E_{1/2})}\nonumber\\
&\leq& |x|_{ E_{1/2}}+\frac{\varepsilon_0}{4}+C|f|_{E_0(J)}\nonumber\\
&\leq&
\frac{\varepsilon_0}{2}+C\left(C\varepsilon_0^2+C\sqrt{T}\right)\nonumber\\
&\leq&\frac{1}{2}\varepsilon_0+C\varepsilon_0^2+C\sqrt{T}\nonumber\\
&\leq&\varepsilon_0,
\end{eqnarray}
provided that $T, \varepsilon$ are both small enough. 
Also by (\ref{E-6}) and (\ref{E-2-3}) we get that
\begin{eqnarray}\label{E-2-8}
|u|_{E_1(J)}&\leq&
\left|\left(e^{-t\triangle_g^2}\right)x\right|_{E_1(J)}+\left|(Kf)(t)\right|_{E_1(J)}\nonumber\\
&\leq&c_1|x|_{E_{1/2}}+C|f|_{E_0(J)}\nonumber\\
&\leq&\varepsilon_0,
\end{eqnarray}
if $\varepsilon, T$ are both small enough. By (\ref{E-2-7}),
(\ref{E-2-8}),  $u\in V_x(J)$.

Now let $x_1, x_2 \in B_{E_{1/2}}(0, \varepsilon)$ be given and $v_1\in V_{x_1}(J),
v_2\in V_{x_2}(J)$. We can get solutions $u_1, u_2$ to 
(\ref{E-2-2}) with initial value $x_1, x_2$ and $f(v_1), f(v_2)$ respectively. It is clear that
\begin{eqnarray}\label{E-2-9}
|u_1-u_2|_{E_1(J)}\leq
\left|e^{-t\triangle_g^2}\left(x_1-x_2\right)\right|_{E_1(J)}+\left|K(f(v_1)-f(v_2))\right|_{E_1(J)}.
\end{eqnarray}
Note
\[f(v_1)-f(v_2)=\triangle_g^2(v_1-v_2)+R_{v_1}-R_{v_2}.\]
For $v_1, v_2$ fixed, we compute similarly as in \eqref{E-2-5}
 that
\begin{equation}\label{E-2-15}
\begin{split}
f(v_1)-f(v_2)=&g^{i\bar j}_{v_2} g^{k\bar l}_{v_2}\p_i\p_{\bar j} \p_k\p_{\bar l} v_2-g^{i\bar j}_{v_1} g^{k\bar l}_{v_1}\p_i\p_{\bar j} \p_k\p_{\bar l} v_1+g^{i\bar j}g^{k\bar l}\p_i\p_{\bar j}\p_k\p_{\bar l} (v_1-v_2 )\\
&+g^{i\bar j}_{v_1} g^{k\bar q}_{v_1} g^{p\bar l}_{v_1}(\p_i g_{p\bar q}+\p_i\p_p\p_{\bar q} v_1)(\p_{\bar j} g_{k\bar l}+\p_j\p_k\p_{\bar l} v_1)\\
&-g^{i\bar j}_{v_2} g^{k\bar q}_{v_2} g^{p\bar l}_{v_2}(\p_i g_{p\bar q}+\p_i\p_p\p_{\bar q} v_2)(\p_{\bar j} g_{k\bar l}+\p_j\p_k\p_{\bar l} v_2)\\
&+g^{i\bar j}\p_{\bar j} g^{k\bar l} \p_i\p_k\p_{\bar l}(v_1-v_2)+g^{i\bar j}\p_{i} g^{k\bar l} \p_{\bar j}\p_k\p_{\bar l}(v_1-v_2)\\&+\left(g^{i\bar j}_{v_2} g^{k\bar l}_{v_2}-g^{i\bar j}_{v_1} g^{k\bar l}_{v_1}\right)\p_i\p_{\bar j} g_{k\bar l}.
\end{split}
\end{equation} 
 
Now we use the similar trick as in \eqref{E-2-10} to deal with \eqref{E-2-15}, then we can get that
\[
\begin{split}
\sqrt{t}|f(v_1)-f(v_2)|_{E_0}\leq &C \sqrt{t}( (|v_1-v_2|_{c^{2, \alpha}}) |v_2|_{c^{4, \alpha}} +|v_2|_{c^{2, \alpha}} |v_1-v_2|_{c^{4, \alpha}})\\
&+C\sqrt{t}\left(|v_1-v_2|_{c^{2, \alpha}}|v_1|^2_{c^{3, \alpha}}+(|v_1|_{c^{3, \alpha}}+|v_2|_{c^{3, \alpha}}) |v_1-v_2|_{c^{3,\alpha}}\right) \\
&+C\sqrt{t}(|v_1-v_2|_{c^{3,\alpha}}+|v_1-v_2|_{c^{2, \alpha}})).
\end{split}
\]
All terms except $|v_1-v_2|_{c^{3, \alpha}}$ in the above inequality are easy to control.  We can estimate
\[
\begin{split}
\sqrt{t}|v_1-v_2|_{c^{3, \alpha}}\leq& \sqrt{t} |v_1-v_2|_{c^{2, \alpha}}^{1/2}|v_1-v_2|_{c^{4, \alpha}}^{1/2}\\
&\leq t^{1/4} (|v_1-v_2|_{c^{2, \alpha}}+\sqrt{t}|v_1-v_2|_{c^{4, \alpha}})
\end{split}
\]

It implies that
\begin{equation}\label{E-2-16}
|f(v_1)-f(v_2)|_{E_0(J)}\leq
C(\varepsilon_0+T^{1/4})|v_1-v_2|_{E_1(J)}+C(\varepsilon_0+{T}^{1/4})|v_1-v_2|_{C(J,
E_{1/2})}.
\end{equation}

We estimate
\begin{eqnarray*}
|v_1-v_2|_{C(J,
E_{1/2})}\leq & |(v_1-v_2)-e^{-t\triangle_g^2}(x_1-x_2)|_{C(J, E_{1/2})}\\
&+|e^{-t\triangle_g^2}(x_1-x_2)|_{C(J, E_{1/2})}.
\end{eqnarray*}
By (\ref{E-7}), 
it follows that we get
\begin{eqnarray}\label{E-2-11}
|v_1-v_2|_{C(J,
E_{1/2})}&\leq& C |(v_1-v_2)-e^{-t\triangle_g^2}(x_1-x_2)|_{E_1(J)}\nonumber\\
&&+|e^{-t\triangle_g^2}(x_1-x_2)|_{C(J, E_{1/2})}\nonumber\\
&\leq& C|v_1-v_2|_{E_1(J)}+C|x_1-x_2|_{E_{1/2}}.
\end{eqnarray}

By (\ref{E-2-9}), \eqref{E-2-16} and (\ref{E-2-11}), we can get that
\begin{eqnarray}\label{E-2-12} |u_1-u_2|_{E_1(J)}&\leq& c_1|x_1-x_2|_{E_{1/2}}+\|K\||f(v_1)-f(v_2)|_{E_0(J)}\nonumber\\
&\leq& C|x_1-x_2|_{E_{1/2}}+C(\varepsilon_0+T^{1/4})|v_1-v_2|_{E_1(J)}\nonumber\\
&\leq& C|x_1-x_2|_{E_{1/2}}+\frac{1}{2}|v_1-v_2|_{E_1(J)},
\end{eqnarray}
provided that $\varepsilon, T$ are both small enough. In
particular taking $x_1=x_2=x$, we obtain from (\ref{E-2-12}) that
\[
|u_1-u_2|_{E_1(J)}\leq \frac{1}{2}|v_1-v_2|_{E_1(J)}
\] for  $v_1, v_2\in V_x(J)$.
The fixed point theorem in Banach space implies that the map $\Pi_x$ has a
unique fixed point
\[
\varphi(\cdot, x)\in V_x(J)\subset C^1_{1/2}(J, E_0)\cap
C_{1/2}(J, E_1)
\] for each $x\in B_{E_{1/2}}(0, \varepsilon)$. By (\ref{E-2-12}), it is
clear that
\[
|\varphi(\cdot, x)-\varphi(\cdot, y)|_{E_1(J)}\leq
C|x-y|_{E_{1/2}}.
\]
Note that we also have
\[
\|\varphi(\cdot, x)-\varphi(\cdot, y)\|_{C([0, T],
E_{1/2})}\leq C|\varphi(\cdot, x)-\varphi(\cdot, y)|_{E_1(J)}.
\]

Lastly, we need to show that the solution of (\ref{1-1}) is
unique. Suppose that $u_1, u_2$ are two solutions of (\ref{1-1}).
Let
\[
T_1=\sup \{t\in [0, T]; u_1(s)=u_2(s)~~\forall~~0\leq s<
t\}.
\]
Since $u_1, u_2\in  C^1_{1/2}(J, E_0)\cap C_{1/2}(J, E_1)$ we
conclude that both belong to the set $V_x(J^{*})$ for $J^{*}=[0,
T^{*}]$ provided that $T^{*}$ small enough, where
\[
V_x(J^{*})=\left\{v\in E_1(J^{*});\; v(0)=x, \|v\|_{C(J^{*},
E_{1/2})}\leq \varepsilon_0\right\}\cap B_{E_1(J^{*})}(0,
\varepsilon_0).
\]
The fixed point theorem in Banach space provides a unique solution in
$V_{x}(J^{*})$ and we conclude that $T_1>0$. Assume that
$T_1<T$. It is clear that $u_1(T_1)=u_2(T_1)=y.$ Let
$v_j(t)=u_j(t+T_1)$, $j=1, 2$ with $t\in J_2=[0, T_2]$ for
some $T_2\in (0, T-T_1]$. Then $v_1, v_2\in C^{1}(J_2,
E_0)\cap C(J_2, E_1)$ and $v_1, v_2$ solve
\begin{eqnarray*}
\frac{\partial \varphi}{\partial t}&=&R_\varphi-\underline{R},\nonumber\\
\varphi(0)&=&y.
\end{eqnarray*}
If $T_2$ is small enough, then $v_1, v_2$ belong to the set
$V_y(J_2)$ and we conclude again that $v_1=v_2$ in the set
$V_x(J_2)$. Therefore $u_1=u_2$ for $t\in [0, T_1+T_2]$,
contradiction.
\end{proof}

Theorem \ref{T-3} is a direct consequence of Theorem \ref{T-1-1} and Theorem 4.1 in \cite{Chenhe}.
\begin{proof} Suppose $g$ is a constant scalar curvature metric on $M$. By Theorem \ref{T-1-1}, there exists $T=T(g), c=c(g), \varepsilon=\varepsilon(g)$ such that for $|\varphi_0|_{c^{2, \alpha}}\leq \varepsilon$, there exists a short time solution $\varphi(t)$ in $[0, T]$.
Also we have estimate
\[
t^{1/2}|\varphi(t)|_{c^{4, \alpha}}\leq c|\varphi_0|_{c^{2, \alpha}}.
\]
By  choosing $\varepsilon$ small enough, we can get that $|\varphi(T)|_{c^{4, \alpha}}\leq T^{-1/2}c|\varphi_0|_{c^{2, \alpha}}$ is still small enough and then we can apply Theorem 4.1 in \cite{Chenhe} to finish the proof.
\end{proof}

\section{Extension for $L^\infty$ metrics}
In this section, we prove that on K\"ahler surfaces one can extend the flow once the metrics are bounded in $L^\infty$, using a blowing up argument.

\begin{theo}\label{T-2}Suppose the Calabi flow exists for time $[0, T)$ on a compact K\"ahler surface $(M, [\omega])$. If there exist  constants $c_1, C_2>0$ such that $c_1\omega_0\leq \omega(t)\leq C_2\omega_0$ for any $t\in [0, T)$, then $T=\infty$ and the curvature is uniformly bounded in $[0, T)$. Also the Calabi flow converges to an extremal metric in Cheeger-Gromov sense when $t\rightarrow \infty$. 
\end{theo}
\begin{proof}
For simplicity, we can assume that the initial metric is $\omega(0)=\omega_0$. 
Suppose the curvature blows up when $t\rightarrow T$, then we can pick up $x_i, t_i$ such that $Q_i=\max_{t\leq t_i}|Rm|=|Rm(x_i, t_i)|\rightarrow \infty.$  Consider the blowing up solution of the Calabi flow
\[
g_i(t)=Q_i g(t_i+tQ_{i}^2)
\]
and the pointed manifold $\{M, x_i, g_i(t)\}$.  
We consider volume ratio for the metric $\omega(t)$ before blowing up. Note that the volume ratio is scaling invariant. Denote $B_p(r, t)$ to be the geodesic ball with respect to $\omega(t)$ at $p$ with radius $r$. We assume that $r\leq r_0$ for some constant $r_0>0$. We can pick up $r_0$ such that for any point $p\in M$, $B_p(r_0, 0)$ is in normal neighborhood with respect to $\omega_0$, and so $B_p(r, 0)$ is diffeomorphic to a Euclidean ball for any $r\leq r_0$. For $r\leq r_0$, we can assume that the volume ratio of $\omega_0$ is bounded away from zero by $\l_0$, 
\[
\frac{Vol_0 (B_p(r, 0))}{r^n}\geq \l_0.
\]
Since $c_1\omega_0\leq \omega(t)\leq C_2\omega_0$, it follows that the volume ratio of $\omega(t)$ is also bounded from below for $r\leq r_0$,
\[
\frac{Vol_t (B_p(r, t))}{r^n}\geq  \frac{\l_0c_1^n}{C_2^n}.
\]

For the sequence $g_i(t)$, the curvature is uniformly bounded and the volume ratio is also bounded from below, then the injectivity radius is also bounded away from zero for $g_i(t)$. Also it is clear  that the Sobolev constant for $\omega(t)$ is uniformly bounded away from zero since $c_1\omega_0\leq \omega(t)\leq C_2 \omega_0$. It follows that  the pointed sequence of the Calabi flow solutions $\{M, x_i, g_i(t)\}$ will converge to a Calabi flow solution $\{M_\infty, x_\infty, g_\infty(t)\}$ by the result in \cite{Chenhe1}.  Moreover,  the limit metric $g_\infty(t)$ is a scalar flat ALE K\"ahler metric on $M_\infty$. Now we show that $M_\infty$ is diffeomorphic to $\mathbb{R}^4$. Roughly speaking, $M_\infty$ is just the limit of the geodesic balls $B_{x_i}(r_0, 0)$ after blowing up. But all these balls are actually diffeomorphic to a Euclidean ball, and it follows that $M_\infty$ itself has to be diffeomorphic to $\R^4$.

Recall that for any $p\in M, r\leq r_0$, the geodesic ball $B_p(r, 0)$ is diffeomorphic to a Euclidean ball. Since $c_1\omega_0\leq \omega(t)\leq C_2 \omega_0$,  then 
 \begin{equation}\label{E-3-1}
 B_p(r/\sqrt{C_2}, t)\subset B_p(r, 0) \subset B_p(r/\sqrt{c_1}, t). 
 \end{equation}
 Let us recall the construction of $(M_\infty, x_\infty, g_\infty)$,  which is the geometric limit in the Cheeger-Gromov sense of a subsequence of $(M, x_i, g_i(0)=Q_i g(t_i))$ for $i\rightarrow \infty$. Without loss of generality, we can assume that $(M, x_i, Q_ig(t_i))$ converges to $(M_\infty, x_\infty, g_\infty)$. 
Denote $g_i=Q_ig(t_i)$. By the construction of $M_\infty$, there exists open subsets $U_i$ of $M_\infty$ for each $i$ such that $x_i\in U_i$, $U_i\subset U_{i+1}$ and $M_\infty=\cup U_i$. Also for each $i$, there exists a diffeomorphsim $\varphi_i: U_i\rightarrow \varphi_i(U_i)\subset \{M, x_i, g_i\}$, $\varphi_i(x_\infty)=x_i$ and $\varphi_i^{*}g_i\rightarrow g_\infty$ on any compact subset $K$ of $M_\infty$. Now for any fixed $R\in (0, \infty)$, consider $B_{x_\infty}(R)\subset M_\infty$. 
For $i$ large enough, we can assume that $B_{x_\infty}(R)\subset U_i$.

 We consider $\varphi_i(B_{x_\infty}(R))\subset \{M, x_i, g_i\}$. Since $\varphi_i$ is almost an isometry when $i$ large enough, we can assume that $\varphi_i(B_{x_\infty}(R))$ contains the geodesic ball centered at $x_i$ with radius $R-1$ and is contained in the geodesic ball centered at $x_i$ with radius $R+1$ with respect to metric $g_i$. Before scaling, we can get \[B_{x_i}(Q_i^{-1}(R-1), t_i)\subset \varphi_i(B_{x_\infty}(R))\subset B_{x_i}(Q_i^{-1}(R+1), t_i).\]
By \eqref{E-3-1}, we can get that
\begin{equation}\label{E-3-2}
 B_{x_i}(Q_i^{-1}(R-1)\sqrt{c_1}/\sqrt{C_2}, t_i)\subset B_{x_i}(Q_i^{-1}(R-1) \sqrt{c_1}, 0)\subset  B_{x_i}(Q_i^{-1}(R-1), t_i).
\end{equation}
It is clear that for $i$ large enough, $B_{x_i}(Q_i^{-1}(R-1) \sqrt{c_1}, 0)$ is diffeomorphic $\R^4$. Now we fix $i=i_R$ large enough. Let $V_R=\varphi_i^{-1}(B_{x_i}(Q_i^{-1}(R-1) \sqrt{c_1}, 0)).$ Then $V_R$ is an open set on $M_\infty$ and it is diffeomorphic to $\R^4$. 
In particular, by \eqref{E-3-2}, 
\begin{equation}\label{E-3-3}
B_{x_\infty} ((R-1)\sqrt{c_1}/\sqrt{C_2}-1) \subset V_R\subset B_{x_\infty} (R).
\end{equation}
By \eqref{E-3-3}, we can pick up an increasing sequence $R_k\rightarrow \infty$ when $k\rightarrow \infty$ 
and correspondingly we can construct $V_{R_k}$ such that $V_{R_k}\subset V_{R_{k+1}}$. It is also clear that
\[
M_\infty=\cup_k V_{R_k}.
\]

Hence $M_\infty$ is diffeomorphic to $\mathbb{R}^4$. But a scalar flat ALE (asymptotically locally Euclidean) K\"ahler metric on $\mathbb{R}^4$ has to be flat \cite{Anderson}. This contradicts that $(M_\infty, g_\infty)$ is non-flat.  This implies that  $T=\infty$ and the curvature has to be bounded for all $[0, \infty)$.  Also it is clear that the flow converges to an extremal metric in Cheeger-Gromov sense.
\end{proof}

{\bf Acknowledgement}: The author would like to thank  Prof. X. Chen for suggesting him to seek the optimal condition for the local solution to the Calabi flow, and  for constant support and encouragement. The author is also grateful to Prof. S. Angenent for insightful discussions about the parabolic equations, and to Prof J. Chen and Prof A. Fraser for constant support.

\end{document}